\documentclass[12pt]{amsart}
\usepackage{amssymb}
\usepackage{amsmath}
\usepackage{graphicx}

\newtheorem{conjecture}{Conjecture}

\newcommand{\figref}[1]{Figure~\ref{#1}}

\begin{document}

\title[Cycles of Nonlinear Dynamical Systems]{Finding, Stabilizing, and Verifying Cycles of Nonlinear Dynamical Systems}

\author{D. Dmitrishin}
\address{Department of Applied Mathematics\\Odessa National Polytechnic University\\Odessa 65044, Ukraine}
\email{dmitrishin@opu.ua}

\author{I.E. Iacob}
\address{Department of Mathematical Sciences\\Georgia Southern University\\Statesboro, GA 30460, USA}
\email{ieiacob@GeorgiaSouthern.edu}

\author{I. Skrinnik}
\address{Department of Applied Mathematics\\Odessa National Polytechnic University\\Odessa 65044, Ukraine}
\email{skrynnyk@opu.ua}

\author{A. Stokolos}
\address{Department of Mathematical Sciences\\Georgia Southern University\\Statesboro, GA 30460, USA}
\email{astokolos@GeorgiaSouthern.edu}




\begin{abstract}
We present a new solution for fundamental problems in nonlinear dynamical systems: finding, verifying, and stabilizing cycles. The solution we propose consists of a new control method based on mixing previous states of the system (or the functions of these states). This approach allows us to locally stabilize and to find a priori unknown cycles of a given length. Our method generalizes and improves on the existing one dimensional space solutions to multi-dimensional space while using the geometric complex functions theory rather than a linear algebra approach. Several numerical examples are considered. All statements and formulas are given in final form. The formulas derivation and reasoning may be found in the cited references. The article focuses on practical applications of methods and algorithms.

\end{abstract}

\maketitle
\date{Oct 7, 2017}

\section{Introduction}
The problem of cycle detection is one of the most fundamental in Mathematics.  The second part of Hilbert's 16th problem asks what can be said about the number and location of limit cycles of a planar polynomial vector field of degree $n$? This problem appears to be one of the most persistent problems in the famous Hilbert's list, second only to the Riemann $\zeta$-function conjecture.

A fundamental tool of dynamics that is often used for analyzing the continuous time system is a reduction of continuous time flow to its Poincar\'e section which is a discrete system. So, an understanding of the discrete systems case is a significant portion of understanding the general situation.

We are developing a new method for detecting high order cycles in discrete autonomous dynamical systems.
Our method is an alternative to what was developed in Physics literature (c.f. \cite{H,PSFB,SD}). The difference consists of using the geometric complex function theory instead of a linear algebra approach developed by physicists. As an improvement we get sharp estimates for the range of cycle multipliers and universal schemes that are more robust and much easier to apply. Some advantages of implementing such type of schemes for problems  in Physics and specific examples may be found in numerous Physics publications, in particular in the ones mentioned above. Another standard field of applications is Biology,  c.f. \cite{PM}.

The rest of the article is organized as follows. In Section~\ref{sec:closed} we define the problem and formally describe our approach in detail. We describe how to find the characteristic polynomials in Section~\ref{sec:poly} and then we define the geometric stability criteria in Section~\ref{sec:geom}. The stability criteria leads to a few optimization problems of which solutions produce the mixing coefficients we use in our method. We discuss the optimization problems in Section~\ref{sec:opti} and their solutions in Section~\ref{sec:coef}. In Section~\ref{sec:nume} we present some numerical simulation results and conclude with Section~\ref{sec:conc}.

\section{Closed loop systems}\label{sec:closed}
Let us consider the discrete dynamical system
\begin{equation} \label{dsc}
x_{n+1}=f(x_n),\qquad  f: A\to A, \; A\subset \mathbb R^m.
\end{equation}
where $A$ is a convex set that is invariant with respect to the function $f$. Let us assume that the system has an unstable  T-cycle  $(x_1^*,...,x_T^*)$. The cycle multipliers $\mu_1,...,\mu_m$ are defined as the zeros of the characteristic polynomial
\begin{equation} \label{cheq}
{\rm det}\left( \mu I-\prod_{j=1}^T   Df (x_j^*) \right)=0.
\end{equation}

In this proposal we restrict ourselves  to considering  multipliers with negative real part, and for convenience let us consider the following two cases:\\

{\bf Case A:} $\{\mu_1,\ldots,\mu_m\}\in\{\mu\in\mathbb R:\mu\in(-\mu^*,1)\}$

{\bf Case B:} $\{\mu_1,\ldots,\mu_m\}\in\{\mu\in\mathbb C:|\mu+R|<R\}\cup \mathbb D,$\\
where $\mathbb D=\{z:|z|<1\}.$

If the cycles are non-stable, which happens when not all multipliers are in the unit disc of the complex plane, then detecting  the cycles might be a difficult problem. In such a case the iterative procedure does not converge, so one has to change the procedure to a more sophisticated one, such as one based on Newton's method.
But even then the instability may still be an issue.

In our work we suggest {\it  changing  the  system rather than the procedure}.  Namely, for the system \eqref{dsc} let us consider an associated closed loop system in the following form
\begin{equation} \label{closed}
x_{n + 1} = \left( {1 - \gamma } \right)f\left(\sum\limits_{j = 1}^N {a_j} x_{n - jT + T} \right) +  \gamma \sum\limits_{j = 1}^N b_j x_{n - jT + 1}
\end{equation}
where $a_1+...+a_N=b_1+...+b_N=1$ and $0\le\gamma<1.$

It is crucial that the system \eqref{closed} preserves  $T$-cycles of the system  \eqref{dsc}.

The first challenge is to find the corresponding characteristic equation. We will tackle this challenge in the next section.

\subsection{Characteristic polynomials}\label{sec:poly}

 The standard approach for finding the characteristic polynomial is  based on the increase of dimensionality to get the quadratic system and then apply linearization. The characteristic polynomial of system \eqref{closed} has coefficients that include  $a_j, b_j$ and the elements of the Jacobi matrices
$Df (x_j^*).$  The expression is very complicated and not practical.

We suggest the  method developed in \cite{DHKS,A} that allows to write the polynomial in very compact and specific form where only the coefficients $a_j, b_j$ and the multipliers $\mu_j$ are used. Using this method we have found \cite{DHKS,A} that the characteristic polynomial of  $T$-cycle can be written in the following elegant form $\lambda^{NTm}f(1/\lambda)$  where
\begin{equation}\label{charpoly}
f(z)=\prod_{j = 1}^m \left(\left[1 - \gamma p(z)\right]^{T}\right. \\- \left.(1 - \gamma )^T \mu _j  z{\left[q(z) \right]}^{\,T} \right),
\end{equation}
and
$$
q\left( z \right) = {a_1} + {a_2}z + \,...\,\, + {a_N}z^{N - 1}
$$
$$
p\left( z \right) = {b_1}\,z + {b_2}\,z^{N - 1} + \,...\,\, + {b_N}z^N,
$$
 The normalization is $q\left( 1 \right) = p\left( 1 \right) = 1$.

\subsection{Geometric Stability Criteria}\label{sec:geom}

The form \eqref{charpoly} above allows us to state the stability criteria which is generalization of a remarkable observation by Alexei Solyanik \cite[p.7]{AS}. Let us consider an auxiliary function
$$
\Phi \left( z \right) = {\left( {1 - \gamma } \right)^T}\,\frac{{z{{\left( {q(z)} \right)}^T}}}{{{{\left( {1 - \gamma \,p(z)} \right)}^T}}}
$$
and the inversion $z^*=1/\bar z.$

A family of the characteristic polynomials of a T-cycle is Schur stable if and only if the following inclusions are valid
\begin{equation}\label{criteria}
\mu _j \in {\left(\bar{\mathbb C} \backslash \Phi (\overline {\mathbb D})\right)^*}, \quad j = 1, \ldots, m.
\end{equation}
Note that $\left(\bar{\mathbb C} \backslash \Phi (\overline {\mathbb D})\right)^*=\overline{\mathbb C} \backslash \Phi^* (\overline {\mathbb D})$  and that in the case of the open loop system $\gamma=0,$ $T=N=1$ we have $\Phi( z) =z,$ and therefore $\left(\bar{\mathbb C} \backslash \Phi (\overline {\mathbb D})\right)^*=\mathbb D.$ Thus we transferred a standard stability criteria from the open loop systems to the closed ones.

\subsection{Optimization Problem}\label{sec:opti}

The stability criteria leads to a few optimization problems, solutions of which produce the required coefficients and  allow us to state the stability criteria in an analytic form.

\subsubsection{Case $\gamma=0$} In this situation
for  case A the following optimization problem is considered
$$
I^{(T)}_N=\sup_{a_j,b_j}\min_{t\in[0,\pi]}\left\{\Re\left(\Phi(e^{it}) \right): \Im\left(\Phi(e^{it}) \right)=0 \right\},
$$
while for case B the following optimization problem is considered
$$
J^{(T)}_N=\sup_{a_j,b_j}\min_{t\in[0,\pi]}\left\{\Re\left(\Phi(e^{it}) \right) \right\}.
$$

Using the above  definitions, the geometric stability criteria can be written analytically as follows:

{The system \eqref{closed} has a stable $T-cycle$ if
 \begin{equation}\label{muopr}
 (\mu^*)| I_N^{(T)}|\le1\quad\mbox{(case A)}
  \end{equation}
  and
  \begin{equation}\label{muopc}
  (R)(2|J_N^{(T)}|)\le1\quad\mbox{(case B)}.
\end{equation}
}

\subsubsection{Case $\gamma\not=0$} The corresponding problems the become to find supremum with respect to possible parameters of the quantities
\begin{equation}\label{AA}
{\rm length}\left\{(\bar{\mathbb C}\backslash\Phi^*(\mathbb D))\cap (-\infty,1)\right\}
\end{equation}
and
\begin{equation}\label{BB}
{\rm area}\left\{(\bar{\mathbb C}\backslash\Phi^*(\mathbb D))\cap \{z:\Re z<1\}\right\}
\end{equation}
If $T=1,2$ then the supremum in the above formulas approaching infinity when $\gamma$ approaches one. Thus, the problems \eqref{AA} and \eqref{BB} make sense only for $T\ge 3.$

Note, that the choice $\gamma=0$ provides the possible choice of the gain and reduces  \eqref{muopr} to \eqref{AA} and \eqref{muopc} to \eqref{BB}. However, in that situation the admissible region will be very narrow in some places, thus a chance to cover a multiplier is more theoretical then practical. One can make it wider and automatically shorter. Then the choice of the polynomial $p(z)$ and $\gamma$ allows to stretch the better region. Thus, we start with the solution of the optimization problem for $\gamma=0$ and then optimize with respect to $\gamma$ and $p(z).$

The next task is finding solutions for the optimization problems above.

\subsection{Coefficients}\label{sec:coef}

The problem of finding $I_N^{(T)}$ and the optimal coefficients was solved for $\gamma=0,T=1,2$ by methods of harmonic analysis in \cite{DKh, DHS, DKKS}. Working extensively on the understanding of the phenomena we came up with the idea of the magnitude of the values $I_N^{(T)}$ and $J_N^{(T)}$ and the polynomials that might be good candidates for the extrema.

\subsubsection{Construction of the polynomials $q\left( z \right)$}

The moving average operation may be treated as a particular kind of low-pass filter, and can be analyzed with the same signal processing techniques  used for low-pass filters, in general. Low-pass filters provide a smoother form of a signal, removing the short-term fluctuations, and leaving a longer-term trend. Thus, the first source for potential solution polynomials can be the set of polynomials that appear in low-pass filters. The most known and important polynomials are the Buterworth polynomials \cite{B}. In our construction we utilize Buterworth type polynomials to define the intermediate polynomials $\eta_N(z).$ Then we apply a Fej\'er type transformation to obtain the desired polynomials.

Let $T$ and $N$ be positive integers, and let $0<\sigma\le \tau\le2.$
We define the set of points
$$
t_j=\frac{\pi(\sigma+T(2j-1))}{\tau+(N-1)T},
\quad j=1,.., \frac {N-2}2 \;\mbox{(N-even)}, \;
\left(\frac{N-1}2\;\mbox{(N-odd)}\right)
$$
and the following generating polynomials:
$$
\eta_N(z)=z(z+1)\prod_{j=1}^{\frac{N-2}2}(z-e^{it_j})(z-e^{-it_j}),\quad\mbox{N-even};
$$
$$
 \eta_N(z)=z\prod_{j=1}^{\frac{N-1}2}(z-e^{it_j})(z-e^{-it_j}),\quad \mbox{N-odd}.
$$
Writing $\eta_N(z)$ in a standard form
$$\eta_N(z)=z\sum_{j=1}^Nc_jz^{j-1}$$
we can define the following three-parameter family of polynomials
\begin{equation}\label{qzT}
q(z,T,\sigma,\tau)=
K\sum_{j=1}^N\left(1-\frac{1+(j-1)T}{2+(N-1)T}\right)c_jz^{j-1},
\end{equation}
where $K$ is a normalization factor that makes $q(1,T,\sigma,\tau)=1.$ \\

\begin{conjecture}\label{con4}
For any $T$ and $N$
$$
q(-1)=-\frac T{2+(N-1)T}\prod_{j=1}^{\frac{N-2}2}\cot^2\frac{t_j}2,\; \mbox{N-even,}
$$
and
$$
q(-1)=-\prod_{j=1}^{\frac{N-1}2}\cot^2\frac{t_j}2,\; \mbox{N-odd}.
$$
Moreover, for $\tau=\sigma$ and for any positive integer $T$
$$
q(-1)\sim N^{-\sigma/T}
$$
\end{conjecture}

\begin{conjecture}\label{con5}
 For any $N$ and $T$ the polynomials $q(z,T,\sigma,\tau)$ are univalent in $\mathbb D$.
\end{conjecture}

The conjectures would provide the justification to the stabilization scheme for real multipliers and $\gamma=0,$ i.e. in case of absence of the linear part in system \eqref{closed} with quantitative estimates of the range for the multipliers.

It is surprising that the addition of a linear part significantly increases an admissible range for the multipliers while also increases the rate of convergence. The next challenge now is to define the polynomial $p(z).$

\subsubsection{Construction of the polynomials $p\left( z \right)$}

For $T = 1$ one can use $p\left( z \right) =z q\left(z\right)$.

For $T = 2$ one can use $\displaystyle p\left( z \right) =1- \frac1{a_1}(1 - z) q\left( z \right)$.

For $T > 2$ it is admissible to use
$$
p\left( z \right) = \frac{2}{{2N - 1}}\left( {z + \, \ldots \, + {z^{N - 1}} + \frac{1}{2}{z^N}} \right).
$$

\begin{conjecture}\label{con6}
 For any $N$ and $T$ there is a choice of $\gamma$ such that the function $\Phi(z)$ is univalent or typically real in $\mathbb D.$ The largest value of $\gamma$ is a point of interest, it produces the widest region for the multipliers.
\end{conjecture}

If $\Phi(z)$ is {\it univalent} or {\it typically real} (typically real means pre-image of a real value is real) then $\Phi(e^{it})$ has only two points of intersection with real axis, namely $\Phi(1)=1$ and $\Phi(-1).$
Thus, for these functions the optimization problem has the estimate
$$
 |\Phi(-1)|\le|q(-1)|^T.
$$
A corollary of the Conjecture \ref{con5} is the following: Choosing $\gamma=0,$ the coefficients $a_j$  produce a closed loop system \eqref{closed} with  stable T-cycle  if
$$
(\mu^*)|q(-1)|^T<1.
$$

A corollary of Conjecture \ref{con6} is the following: The coefficients $a_j$ and $b_j$ produce a closed loop system \eqref{closed} with  stable T-cycle  if
$$
(\mu^*)|\Phi(-1)|<1.
$$
The above inequalities can be tested numerically and the proposed coefficients can be used in the closed loop system  \eqref{closed} to fulfill the main goal of this work -- to numerically detect cycles of high order. This is addressed in Section \ref{sec:nume}.

Case B has been less analyzed. We state the following main conjecture associated with Case B:

 \begin{conjecture}\label{con7}
 For any $N$ and $T$ the polynomials $q(z)$ with $\sigma=\tau=1$ give the solutions to the optimization problem $J_N^{(T)}.$
\end{conjecture}

The numerical testing results seem to indicate that  Conjecture~\ref{con7} is likely  valid.

\section{Numerical simulations}\label{sec:nume}

We performed numerous numerical simulations, of which results, in our opinion, are an important part of this work. Regardless of the theoretical justification, one can apply the methods developed here to detect cycles. In the sequel we list a number of maps  and cycles detected using our method.

The first example is the H\'enon map. In 2016 in the paper \cite{LZK} cycles of lengths 1, 2, 4, and 6 for the H\'enon map were detected. Using our method, cycles of lengths 11 and 28 for the H\'enon map are  detected and presented  below. Note that 11 is a prime integer, and detecting the  cycles of prime length is a much more subtle issue according to the celebrated Sharkovsky theorem. Thus, it is not a coincidence that no cycles of the length 3 and 5 were mentioned in \cite{LZK}.

\subsection{H\'enon map, n=1,...,1200}
\begin{figure}[ht]
\centerline{\includegraphics[scale=0.35]{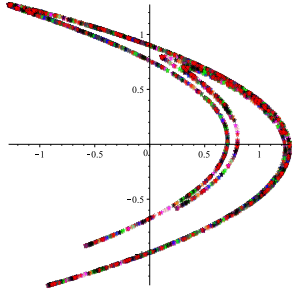}}
\caption{The H\'enon map}\label{fig:henon}
\end{figure}

\begin{figure}[ht]
{\includegraphics[scale=0.25]{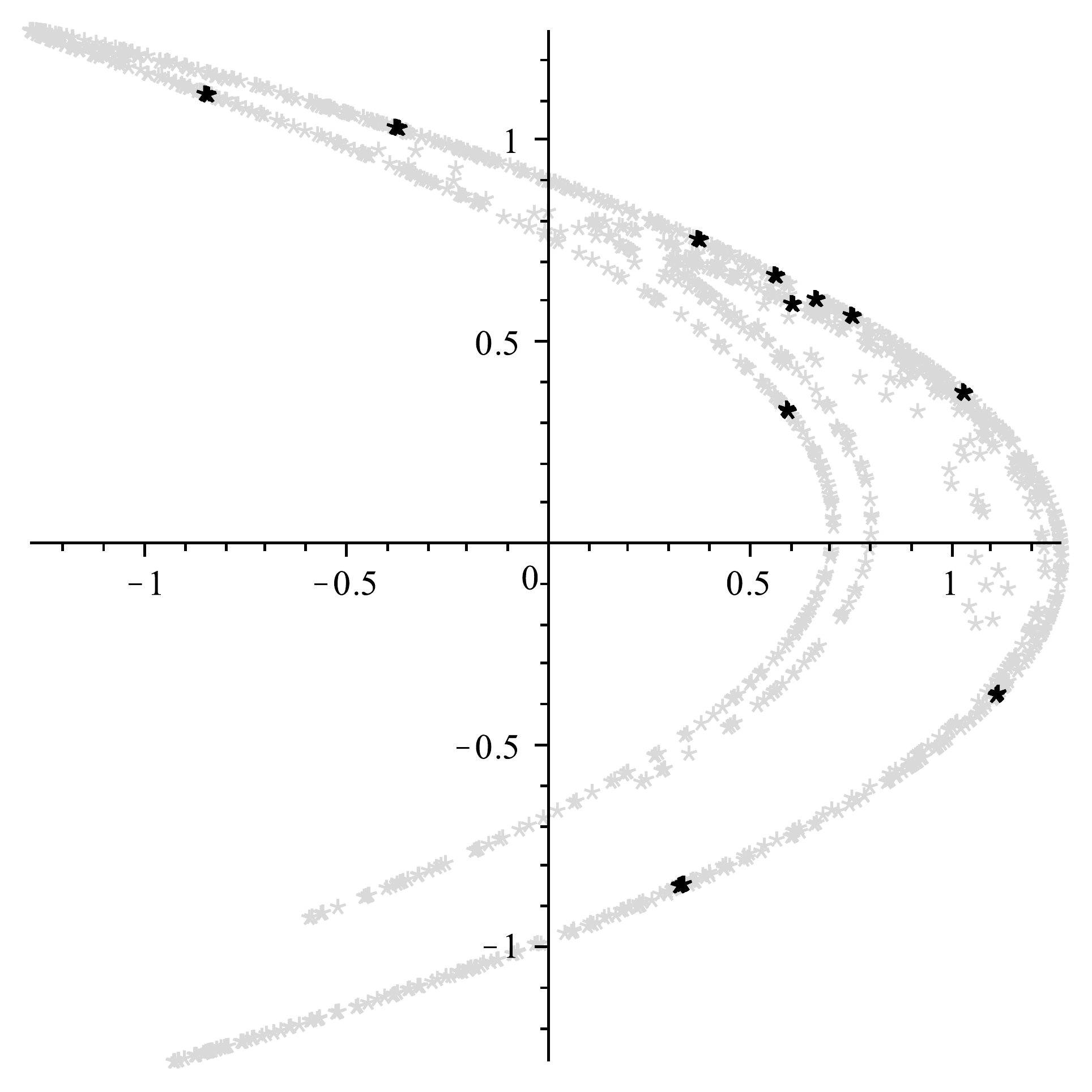} \hspace{.4cm}
\includegraphics[scale=0.25]{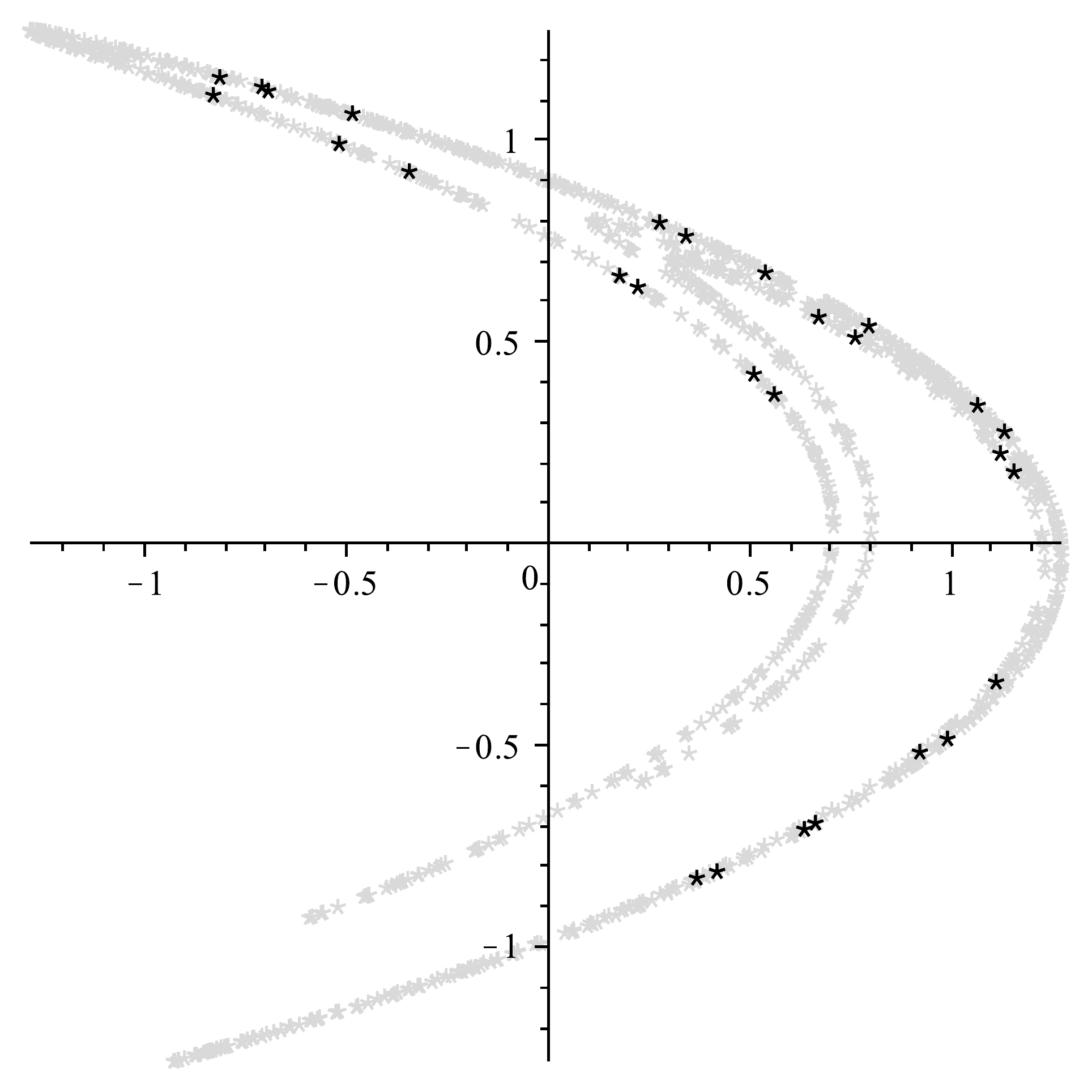}}
\caption{The stabilized H\'enon map. Left: $T = 11, N=10, n=9700,\dots,9900$. Right: $T = 28, N=40, n=18900,\dots,18990$}\label{fig:henont11-28}
\end{figure}

The H\'enon map is described by the system:
$$
\begin{cases}
x_{n+1}=&1-1.4x_n^2+y_n\\
y_{n+1}=&0.3x_n
\end{cases}
$$
\figref{fig:henon} shows the H\'enon chaos and \figref{fig:henont11-28} shows the stabilized H\'enon map for $T = 11$ and $T = 28$, respectively.

\subsection{Elhadj-Sprott map}
\begin{figure}[ht]
{\includegraphics[scale=0.25]{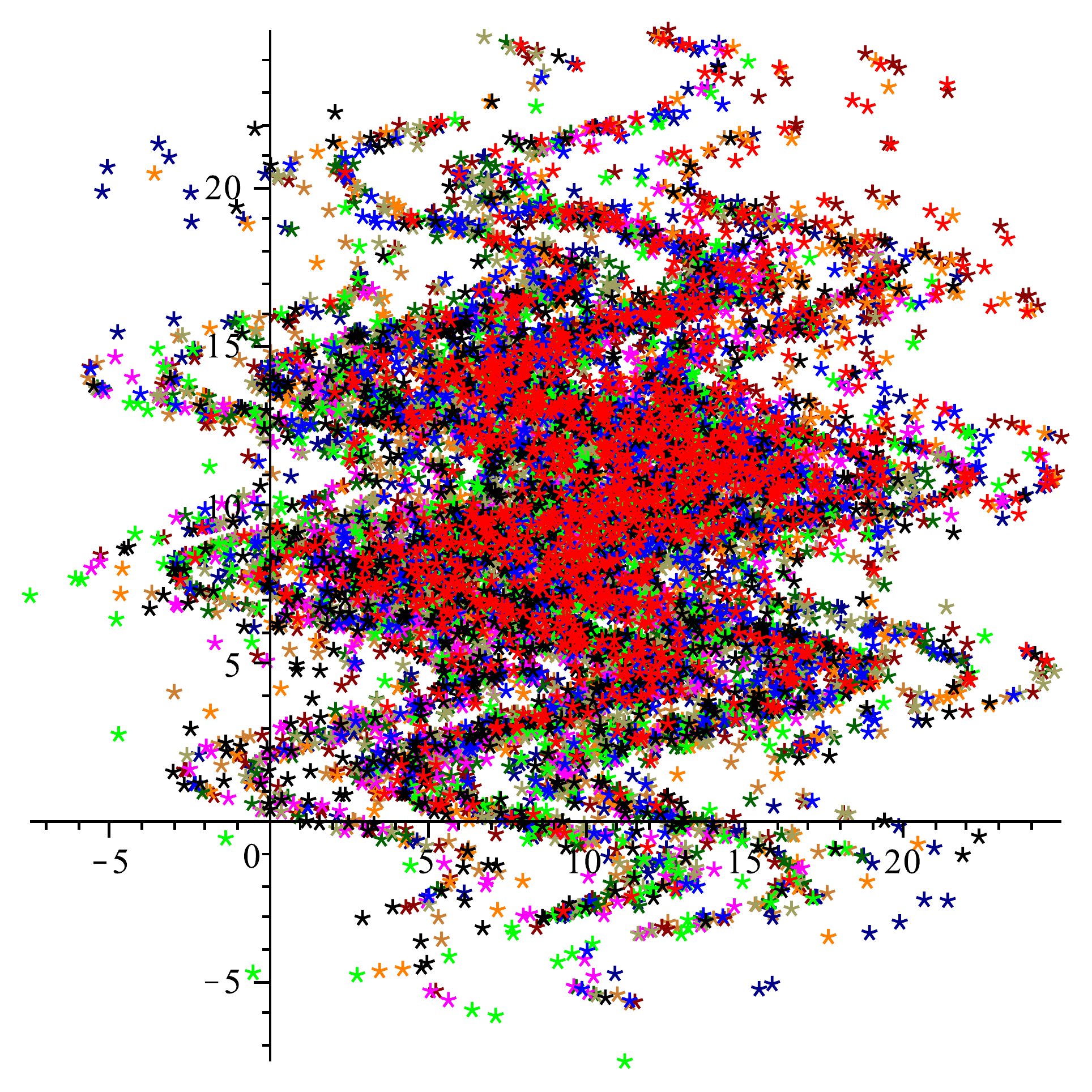} \hspace{.8cm}
\includegraphics[scale=0.25]{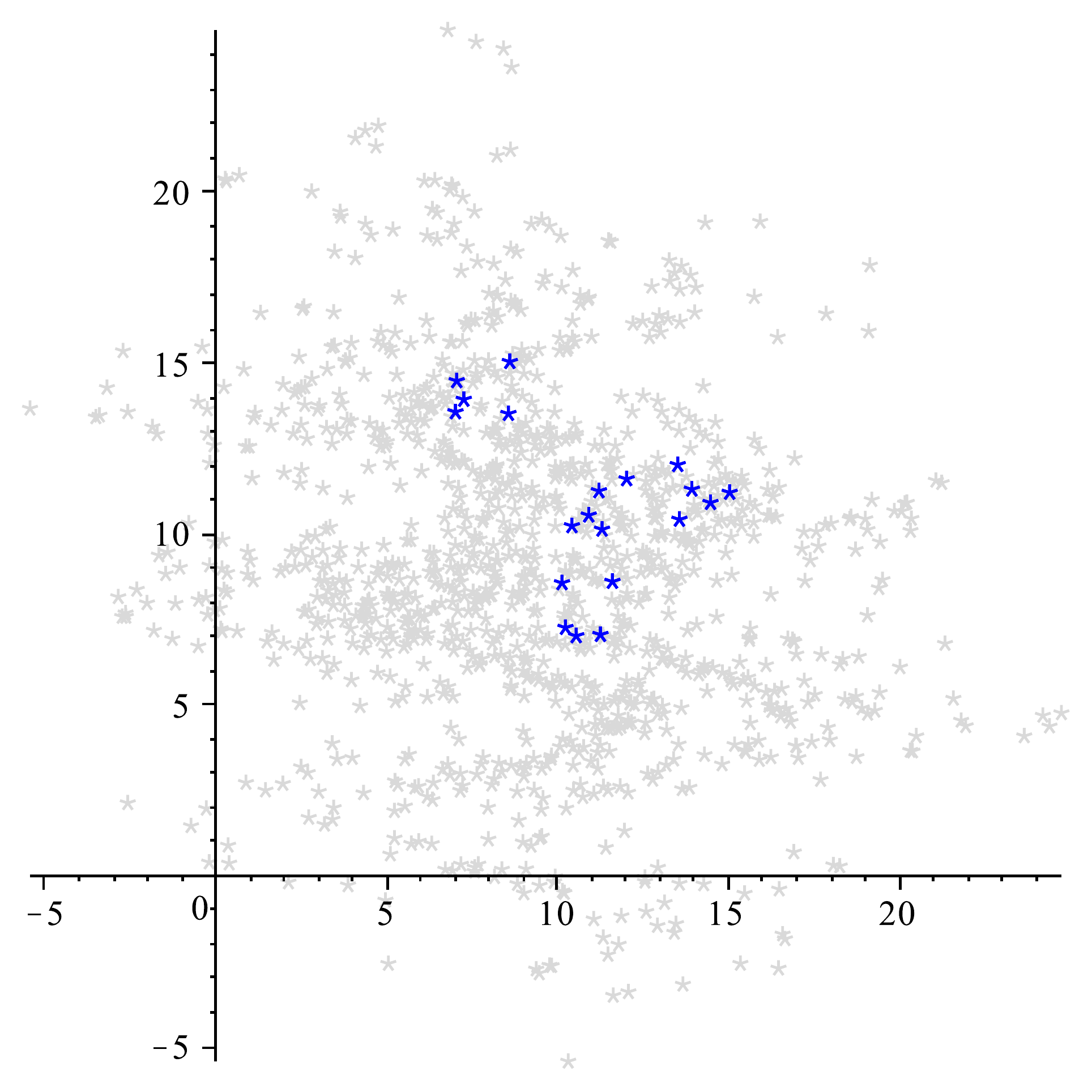}}

\caption{Elhadj-Sprott chaos (left) and its stabilized map for $T=20, N=38, n=23900,\dots,23990$ (right)}\label{fig:sprottt20}
\end{figure}

The Elhadj-Sprott map is described by the system:
$$
\begin{cases}
x_{n+1}=1-4\sin(x_n)+0.9y_n \\
y_{n+1}=x_n
\end{cases}
$$
The Elhdj-Sprott chaos map and its corresponding stabilized version are shown in \figref{fig:sprottt20}, left and right, respectively.

\subsection{Ikeda map}
\begin{figure}
{\includegraphics[scale=0.4]{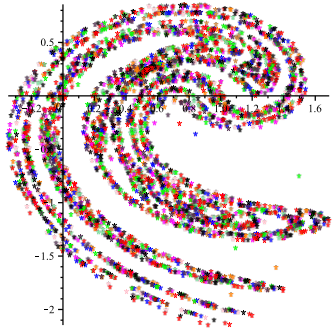}}
\caption{The 2D Ikeda map}\label{fig:ikeda}
\end{figure}

The Ikeda map is described by the system:
$$
\begin{cases}
x_{n+1}=&1+0.9\left(x_n\cos\left(0.4-\frac6{1+x_n^2+y_n^2}\right) - y_n \sin\left(0.4-\frac6{1+x_n^2+y_n^2}\right) \right),
\\
y_{n+1}=& 0.9\left(x_n\sin\left(0.4-\frac6{1+x_n^2+y_n^2}\right) + y_n \cos\left(0.4-\frac6{1+x_n^2+y_n^2}\right) \right),
\end{cases}
$$
\figref{fig:ikeda} shows the Ikeda chaos.

\begin{figure}[ht]
{\includegraphics[scale=0.25]{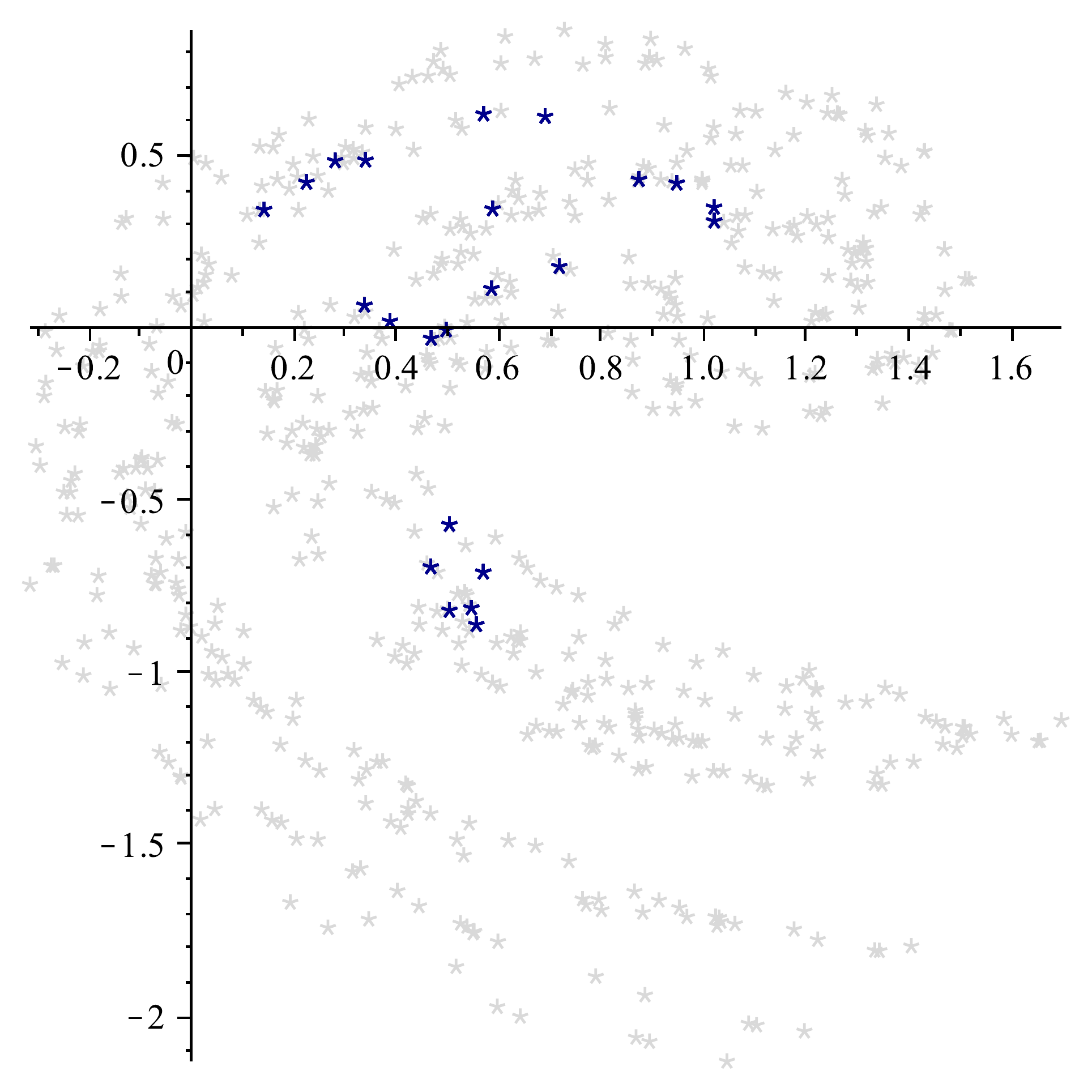}}

\caption{The 2D Ikeda map for $T=23, N=36, n=16700,\dots,16990$}\label{fig:ikedat23}
\end{figure}

\bigskip
The 23 cycle of the Ikeda map is shown in \figref{fig:ikedat23}.

The 5 closest values to the 23-periodic point are listed below:
\bigskip

       1, 0.28041732592998354255, 0.48338110785346721899

       2, 0.28041730651756333896, 0.48338109279633566551

       3, 0.28041728677011750510, 0.48338107747923293207

       4, 0.28041726684498716885, 0.48338106202421746957

       5, 0.28041724714287064958, 0.48338104674212881602

\bigskip

Note that even checking that $x= 0.280417$ and $y= 0.483381$ are truncations of the decimals of the 23-periodic point is a challenge.

\subsection{Lozi map}
\begin{figure}
{\includegraphics[scale=0.25]{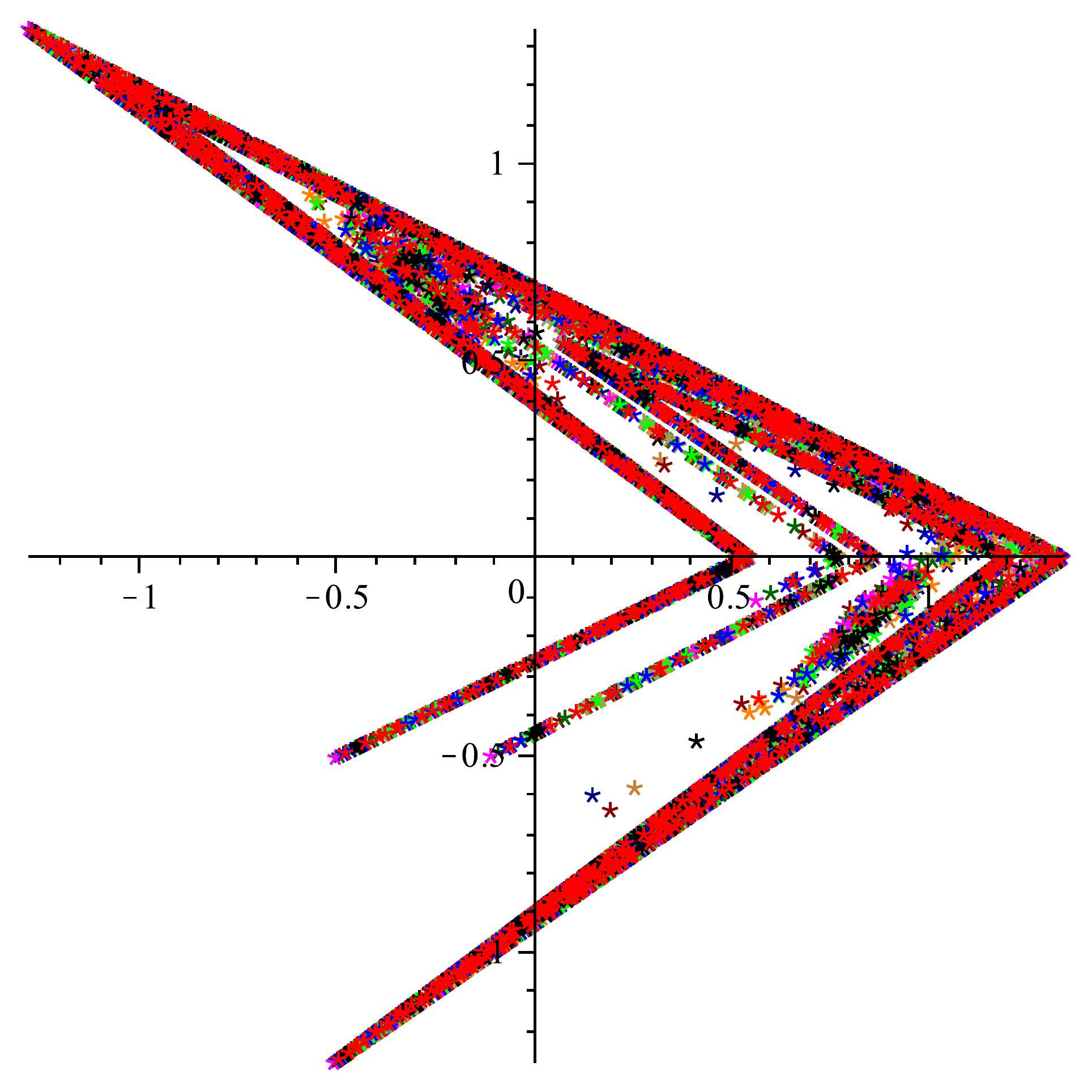} \hspace{.4cm}
\includegraphics[scale=0.25]{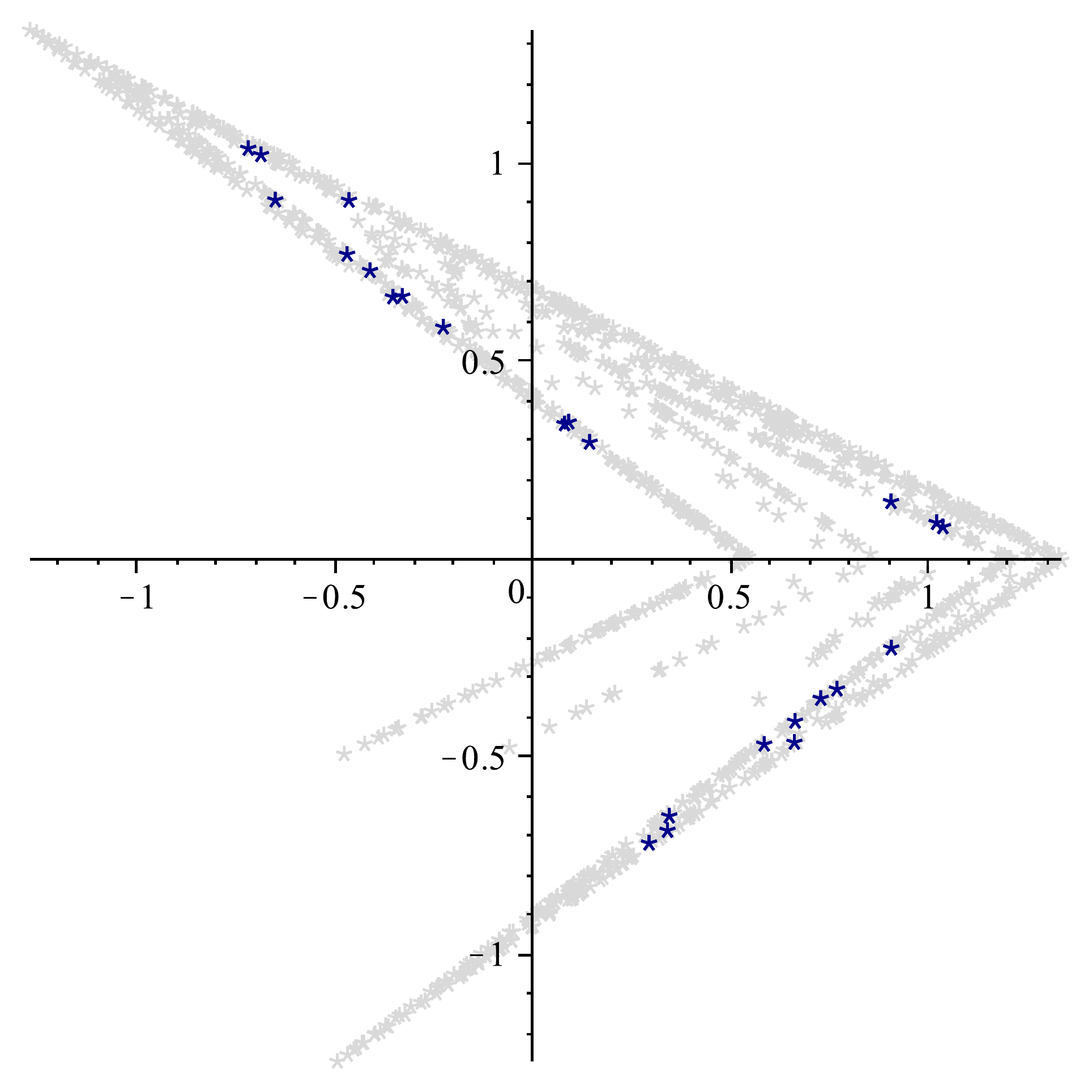}}
\caption{Lozi chaos (left) and the stabilized 24 cycle of the Lozi map for $T=24, N=20, n=5200,\dots,6050$ (right)}\label{fig:lozi}
\end{figure}

The Lozi map is described by the system:
$$
\begin{cases}
x_{n+1}=1-1.7|x_n|+0.5y_n \\
y_{n+1}=x_n
\end{cases}
$$
The Lozi chaos map and its corresponding stabilized version are shown in \figref{fig:lozi}, left and right, respectively.

\subsection{Holmes cubic map}
\begin{figure}[ht]
{\includegraphics[scale=0.3]{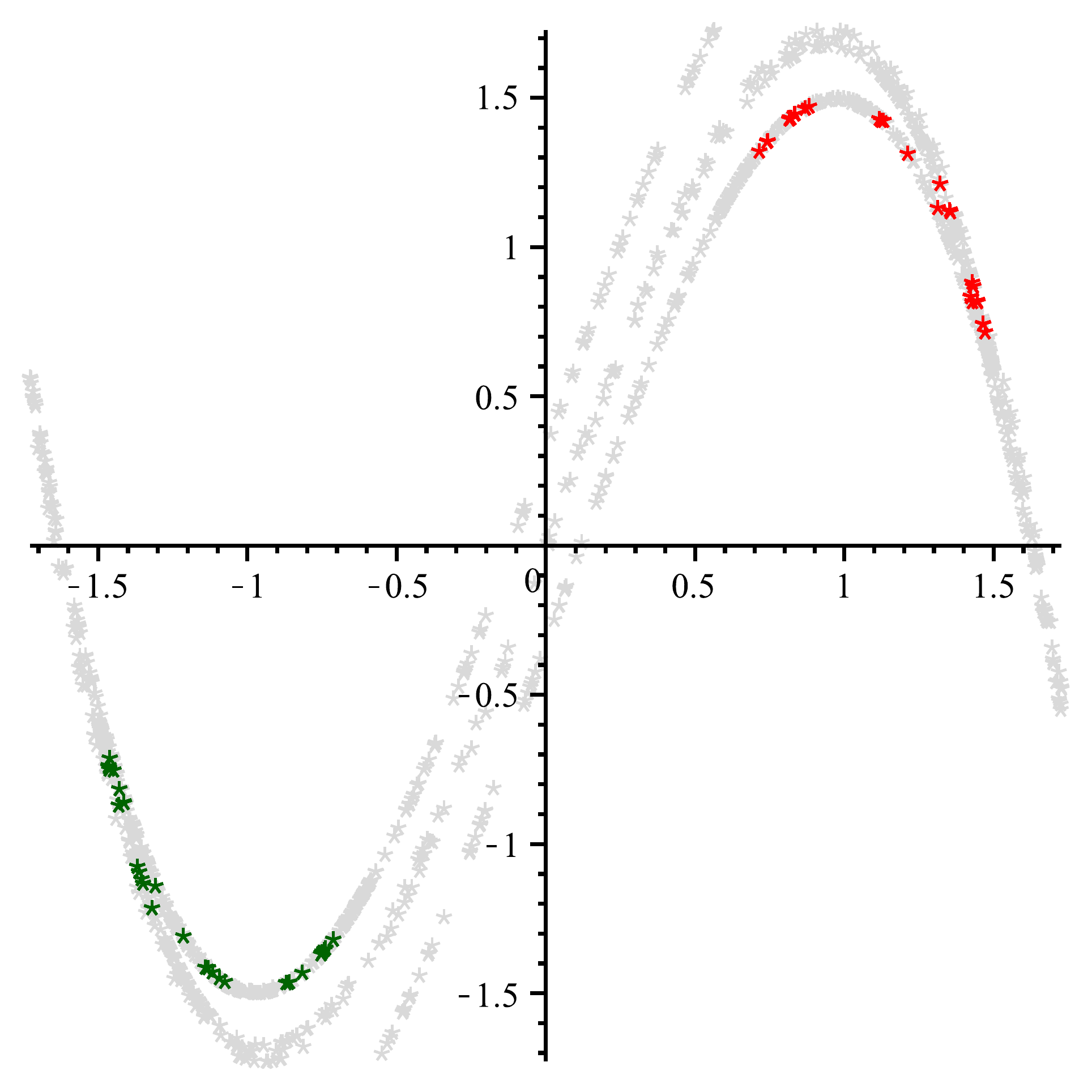}}
{\includegraphics[scale=0.3]{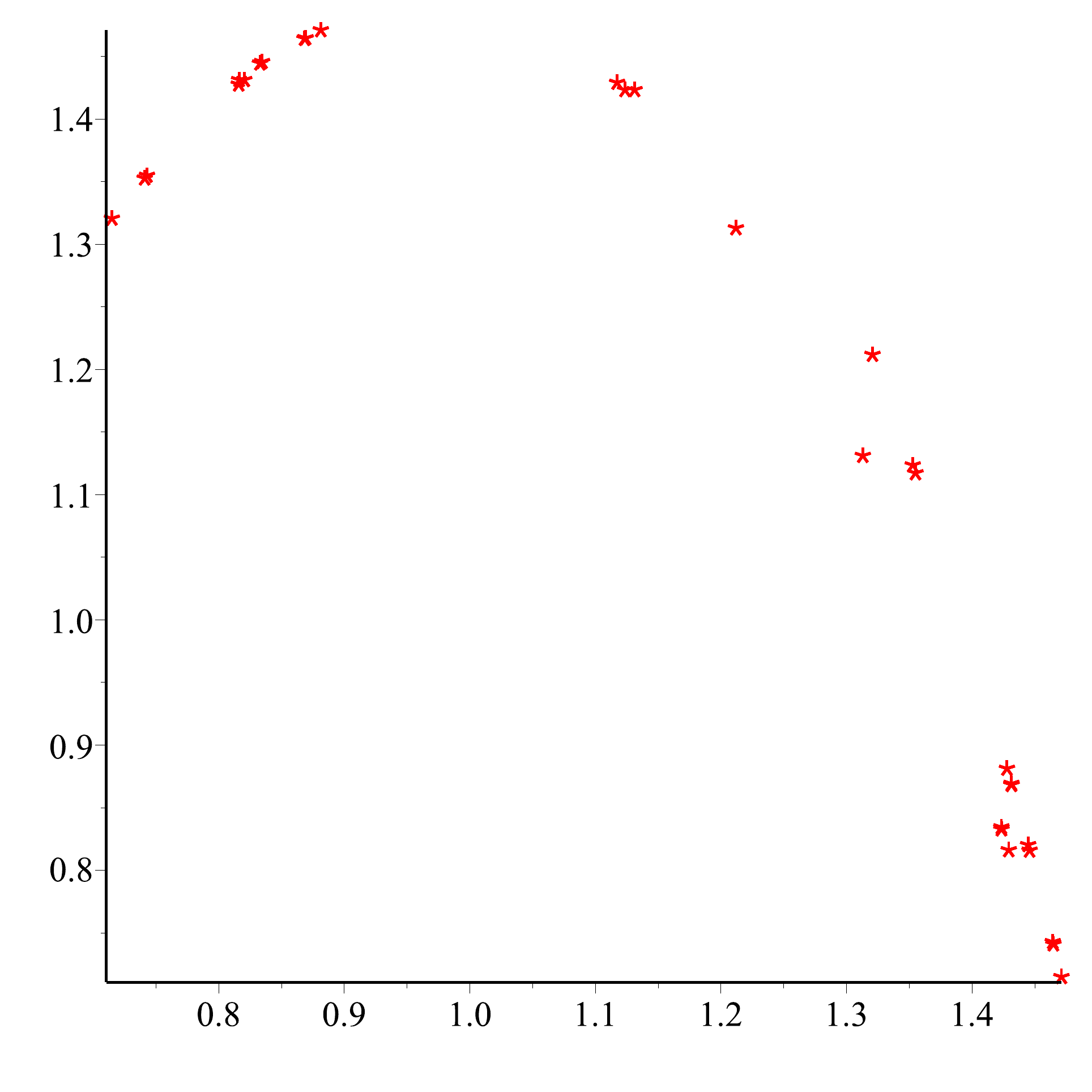}}
\caption{Holmes chaos (left) and the stabilized 24 cycle of the Holmes map for $T=30, N=34, n=5200,\dots,6050$ (right)}\label{fig:holmes}
\end{figure}

Finally, the up to date record is the Holmes cubic map described by the system:
$$
\begin{cases}
x_{n+1}=&y_n \\
y_{n+1}=&1-0.2x_n+2.77y_n - y_n^3
\end{cases}
$$
The Holmes chaos map and its corresponding stabilized version are shown in \figref{fig:holmes}, left and right, respectively.

\subsection{Numerical  difficulties}

In this section we describe some of the numerical difficulties that are likely to be encountered.
To overcome these difficulties one needs to carefully implement the  algorithmic methods. For example, the problem to determine 30-cycle in the Holmes cubic map is equivalent to finding  a solution to a polynomial equation of degree $3^{30}$. Even to verify that a  given number from the cycle has correct digits can be a challenge.

Since our approach is multi-parametric,  an optimization over the parameters can be performed,  as shown in the diagrams below, which leads to significant computational performance improvement of the method.

Furthermore, an interesting phenomena has been observed: the increase of the depth of the used prehistory does not necessarily improve the situation. On the contrary, it   definitely makes things worse when parameter $N$ is large enough. That  is the motivation to look for new schemes that are based only on a few elements from prehistory.

In what follows we list several specific challenges.

\begin{figure}[ht]
{\includegraphics[scale=0.45]{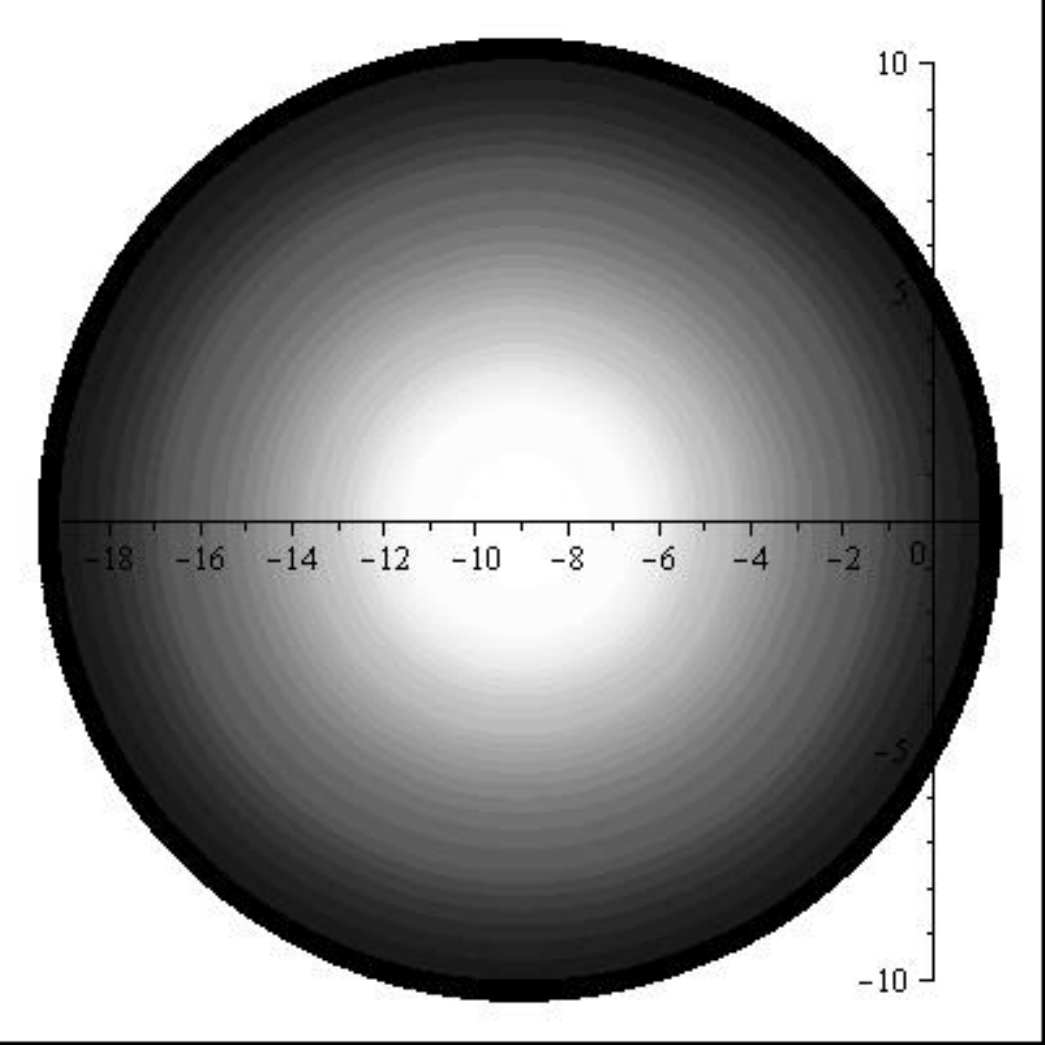}}
\caption{The set $\left(\bar{\mathbb C} \backslash \Phi (\overline {\mathbb D})\right)^*$: distances to the boundary of the unit disk (darker means closer)}\label{fig:phifunc}
\end{figure}

First, the rate of convergence depends on the multipliers distance to the boundary of the region of convergence. In the simplest case scenario  $T=N=1$ the function $\Phi(z)$ is $\displaystyle\Phi(z)=(1-\gamma)\frac z{1-\gamma z}$. \figref{fig:phifunc} displays the set
$\left(\bar{\mathbb C} \backslash \Phi (\overline {\mathbb D})\right)^*$  with $\gamma=0.9$. Different shades  indicate the multiplier distance to the boundary of the unit disc $\mathbb D$ in the  {\it closed loop} system \eqref{closed}. More specifically, the darker the region is the closer the multiplier of the closed loop system is to the boundary $\partial \mathbb D$ , and therefore the {\it  convergence is slower }.

In this case,  as $\gamma$  approaches one, the set $\left(\bar{\mathbb C} \backslash \Phi (\overline {\mathbb D})\right)^*$ can cover any given multiplier with negative real part. However, the white region of good convergence will be centered at $\displaystyle\left(\frac\gamma{\gamma-1},0\right)$, therefore, if we have  small and large in absolute value multipliers, then unavoidably, one of the multipliers will be in the dark zone. Hence, the rate of convergence will be slow regardless of the choice of $\gamma$. This is another manifestation of the stiffness effect in the numerical computations.

\begin{figure}[ht]
{\includegraphics[scale=0.44]{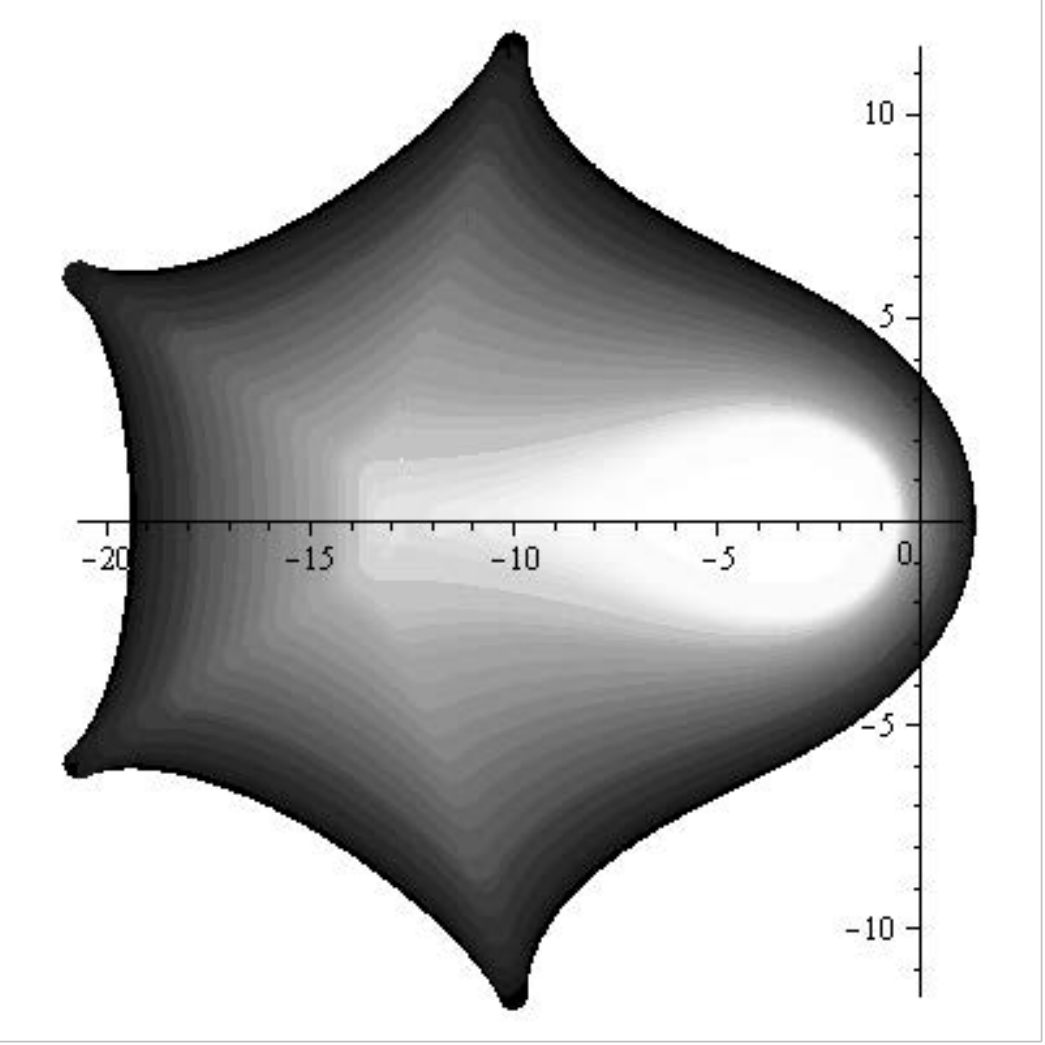}\hspace{1cm}
\includegraphics[scale=0.44]{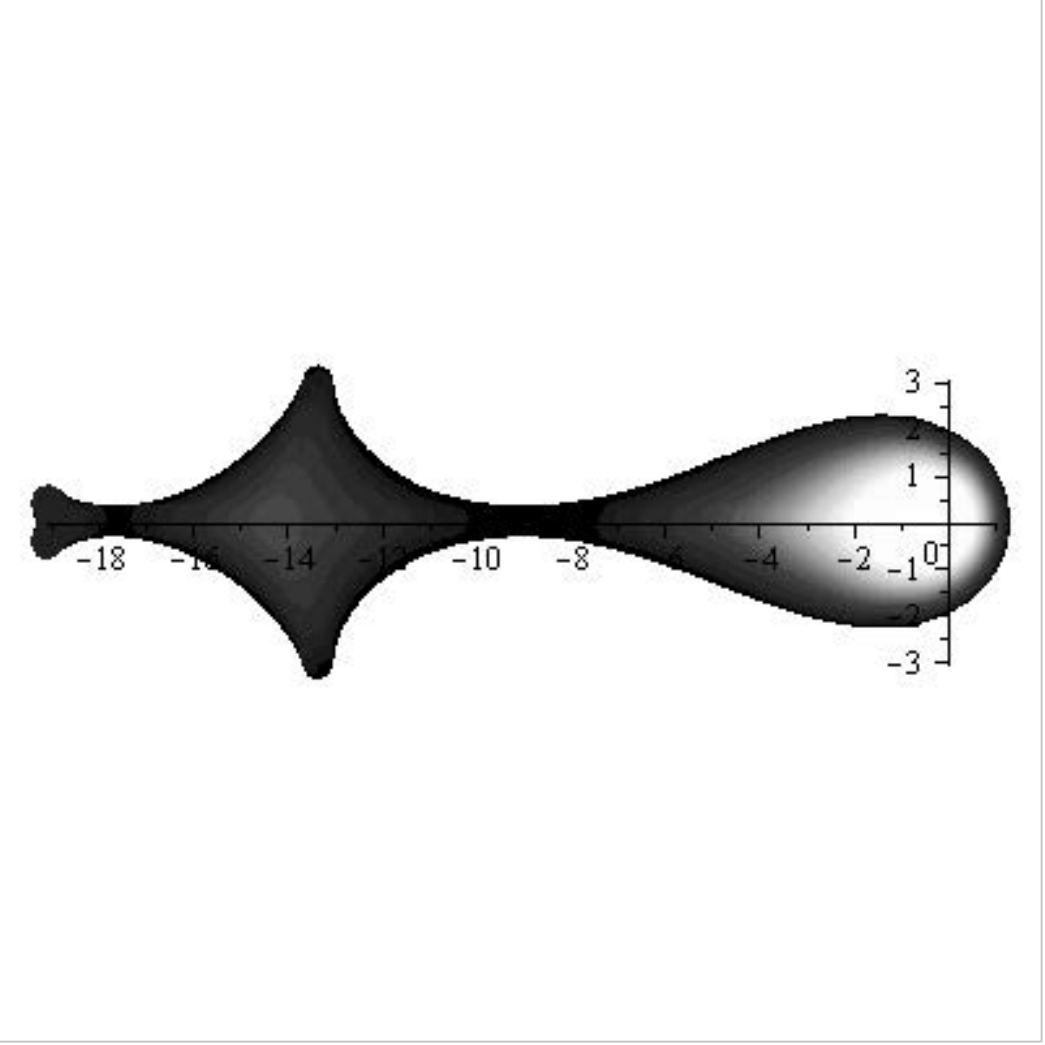}}
\centerline{$\gamma=0.7, \sigma=1$\hspace{2.2cm}$\gamma=0.254, \sigma=2$ }

{\includegraphics[scale=0.44]{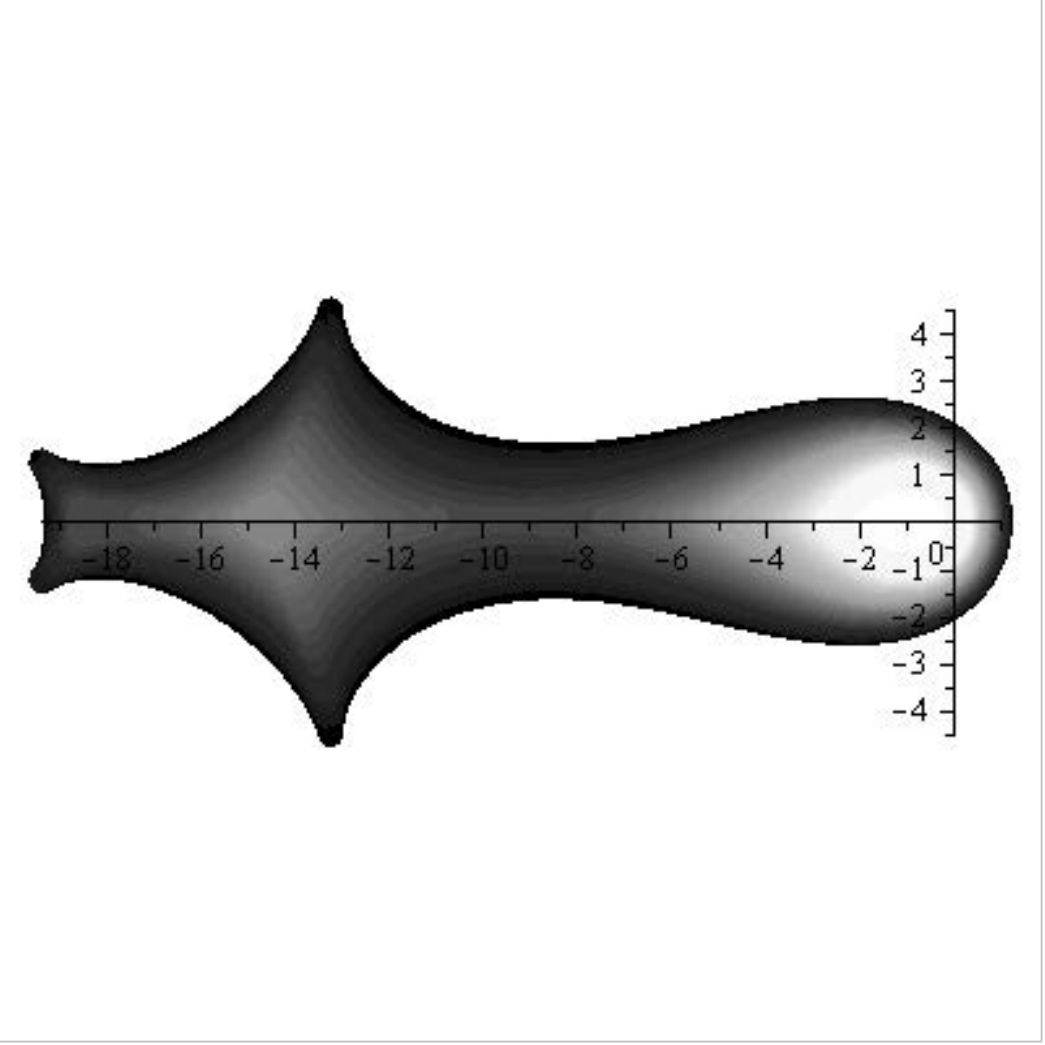}\hspace{1cm}
\includegraphics[scale=0.44]{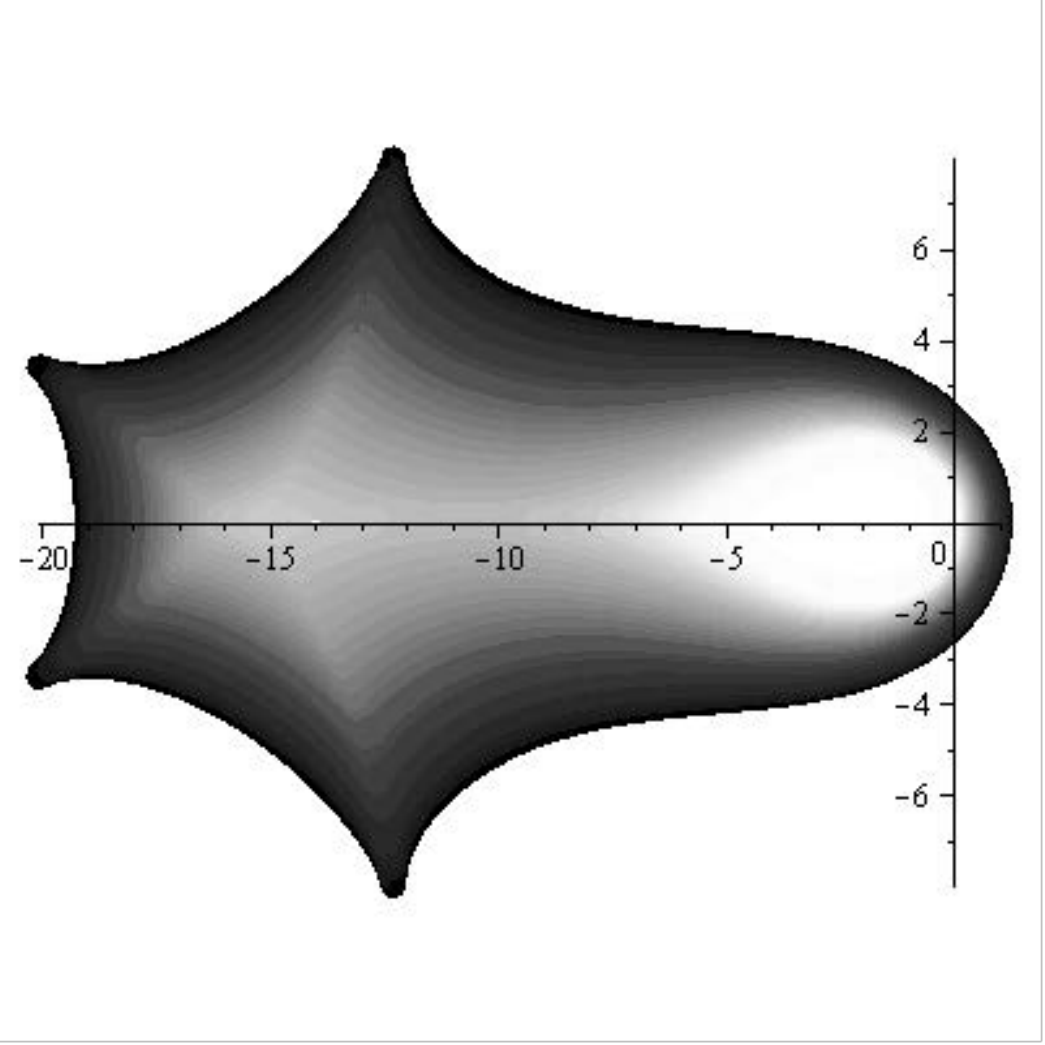}}
\centerline{$\gamma=0.368, \sigma=1.8$\hspace{2.2cm}$\gamma=0.557, \sigma=1.4$ }
\caption{The effect of changes in parameters $\gamma$ and $\sigma$.}
\label{fig:effectchggamma}
\end{figure}

For Figures~\ref{fig:effectchggamma} we have $N=5$ and $\displaystyle\Phi(z)=(1-\gamma)\frac {zq(z)}{1-\gamma zq(z)}.$ The figures show the effect of  changes in the parameters $\sigma$ and $\gamma$. Letting $\sigma$  approach zero make the regions shorter along real axis, and taller along the imaginary axis, therefore the white spot is wider. The same case is if $\gamma$ approaches zero. Recall that  Conjecture \ref{con4} states that the width of the region is about $N^{\sigma}$,  which leads to the optimization problem involving  parameters $N, \sigma$, and $\gamma$.

It is shown in \cite{DFS} that in the case when $T=1$ the largest value $\mu^*$  that allows to fit the multipliers of the closed-loop system into the disc of radius $0< \rho <1$  is asymptotically  about $\displaystyle\frac{4\rho}{(1-\rho)^2}$, and the largest value of $R$ that allows to fit the multipliers of the closed-loop system into the disc of radius $0<\rho<1$  is about $\displaystyle\frac{\rho}{1-\rho}$. In particular, if $\mu<-3600$ then there is no way to fit the multipliers of the closed loop system in a disc of radius 0.9.

Furthermore, the size and shape of the regions of convergence depend on $N$ and they are not nested. Therefore, increasing $N$ does not guarantee  the improvement of convergence. Moreover, in some cases one can observe the change of behavior from stable to chaotic when $N$ increases.

Second, even when solutions are obtained, it is difficult to verify them. It does not help  to substitute the obtained solution  as an initial value to the initial system because of instability. For the same reason it is not recommended  to substitute them in the $T$-iterated system.

Third, the number of iterations is pretty high, therefore the rounding error is a serious issue, especially, having in mind that the coefficients are non-negative numbers  and that their sum has to be one. For example, if we are looking for the 50-cycle, half of the coefficients will have small values, very close to zero.

\section{Conclusion}\label{sec:conc}
In this article we discuss the problem of finding and verification periodic non-stable orbits in non-linear systems in discrete time. As opposed to existent solutions of this problem based on algebraic methods (reduced to solving a system of non-linear equations) we suggest dynamic system approach. Namely, we construct an auxiliary dynamic system for which the periodic orbits coincide with the ones of the original system. However, the periodic orbits of the new system became locally asymptotically stable.

The advantages of a dynamic system approach can be easily illustrated for the example of a simple logistic equation $x_{n+1}=\mu x_n(1-x_n)$ where $\mu$ is slightly smaller than 4. Say, we want to find a cycle of length 20. The algebraic approach leads to the problem of finding the real roots of polynomials of degree $2^{20}$ in the interval $[0,1].$ Let us assume that we want to find the periodic orbit with accuracy $10^{-10}$. However, the roots of the polynomial equations may be closer to each other than $10^{-10}$. A natural question then arises: how to check whether a given root corresponds to a given orbit. The algebraic approach would be a poor choice for making such a verification. If we use the dynamic approach then the obtained points can be used as initial values. If these approximate values correspond to the cycle, then the initial values are in the basin of attraction of the cycle. Our suggested procedure allows us to verify whether that is the case. We describe some classic model equations as practical examples.

Ren\`e Lozi, a well-known expert in non-linear dynamics, posted the following question  \cite{L}: ``{\it Can we trust in numerical computations of chaotic solutions of dynamical systems?}". He concluded:
{\it
``We have shown, in the limited extend of this article, on few but well known examples, that it is very difficult to trust in numerical solution of chaotic dynamical dissipative systems. In some cases one can even proof that it is never possible to obtain reliable results."}

The methodology developed in our presentation allows performing numerical simulations with confidence.


\end{document}